\documentclass{article}
\usepackage[T1]{fontenc}
\usepackage{graphicx}
\usepackage{amsfonts}
\usepackage{amsmath}
\usepackage{url}
\usepackage{hyperref}
\usepackage{cite}

\begin{document}
\title{Numerical Aspects of Hyperbolic Geometry}
%
%
\author{
Dorota Celi\'nska-Kopczy\'nska and Eryk Kopczy\'nski
\\ \small{Institute of Informatics, University of Warsaw}
}

\maketitle              
\begin{abstract}
  Hyperbolic geometry has recently found applications
  in social networks, machine learning and computational biology.
  With the increasing popularity, questions about the best representations of hyperbolic spaces arise, 
  as each representation comes with some numerical instability. This paper compares various 2D and 3D
  hyperbolic geometry representations. To this end, we conduct an extensive simulational scheme based
  on six tests of numerical precision errors. Our comparisons include the most popular models and
  less-known mixed and reduced representations.
  According to our results, polar representation wins, although the halfplane invariant is also very successful.
  We complete the comparison with a brief discussion of the non-numerical advantages of various representations.
  
\end{abstract}
\def\bbR{\mathbb{R}}
\def\bbC{\mathbb{C}}
\def\bbU{\mathbb{U}}
\def\bbH{\mathbb{H}}
\def\bbE{\mathbb{E}}
\def\bbS{\mathbb{S}}
\def\bbV{\mathbb{V}}
\def\bbP{\mathbb{P}}
\def\acos{\mathrm{acos}}
\def\acosh{\mathrm{acosh}}
\def\midp{{\mathrm{mid}}}

\section{Introduction}

Hyperbolic geometry has recently gained interest in many fields. Notable examples include
the hyperbolic random graph models of hierarchical
structures \cite{munzner,lampingrao}, social network analysis \cite{bogu_internet},
hyperbolic embeddings used in machine learning \cite{nickel}, as well as 
visualizations and video games \cite{hyperrogue,hyperbolica,dziadICCS}.

One crucial aspect of hyperbolic geometry is its \emph{tree-like structure},
as shown in Figure 1. Each edge of these trees has the same length,
making them grow exponentially. This tree-likeness property is crucial in 
the modeling hierarchical data \cite{bogu_internet,nickel} and game design
\cite{hyperrogue}. However, this property comes with a severe numerical cost
\cite{tobias_alenex,dhrgex,mltiles,reptradeoff}. Since the circumference of
a hyperbolic circle of radius $r$ is exponential in $r$, any representation
based on a fixed number $b$ of bits will not be able to distinguish between
points in a circle of radius $r=\Theta(b)$, even if the pairwise distances between
these points are large.

Different communities use different representations of hyperbolic spaces. 
The newcomers and some experts often use the Poincar\'e model, most popularly
used in the introductions to hyperbolic geometry. However, in visualizations
\cite{hyperrogue,hyperbolicvr,munzner} the Minkowski hyperboloid model seems to be
commonly used for the internal representation (and converted to Poincar\'e
for visualization purposes), and in the social network research, native
polar coordinates are popular \cite{bogu_internet,friedrich2023computing}.
The users motivate their choice
of representation by factors such as ease of use, generalizability,
and numerical stability. The numerical aspect needs further study; for example,
the hyperboloid model and the Poincar\'e
disk model may be better numerically depending on the computation at hand.
At the moment, we know one study comparing the hyperboloid model and
the Poincar\'e half-plane model \cite{achilles}; however, this paper compares
only two representations of isometries, and has been written in 2002, before
the surge of interest in hyperbolic geometry.

In this paper, we compare
a large number of representations of hyperbolic geometry. We primarily focus on the numerical issues. It is worth to
note that, due to the exponential growth, any representation based on a fixed
number of bits will introduce numerical errors, a commonly used solution to
this \cite{wpigroups,hyperrogue,mltiles,dhrgex,gentes,rtvizfinal} is to use combinatorially
generated tessellations, and represent points and isometries by a pair $(t,h)$,
where $t$ is one of the tiles of the tessellation, and $h$ is coordinates
relative to tile $t$. Our research takes this into account.

\begin{figure}
\begin{center}
\includegraphics[width=.19\textwidth]{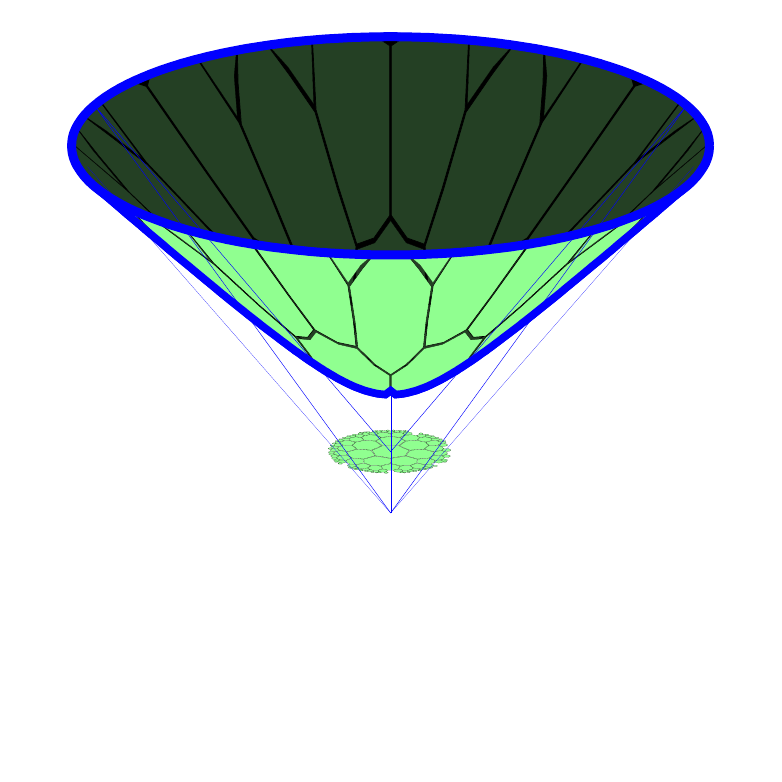}%
\includegraphics[width=.19\textwidth]{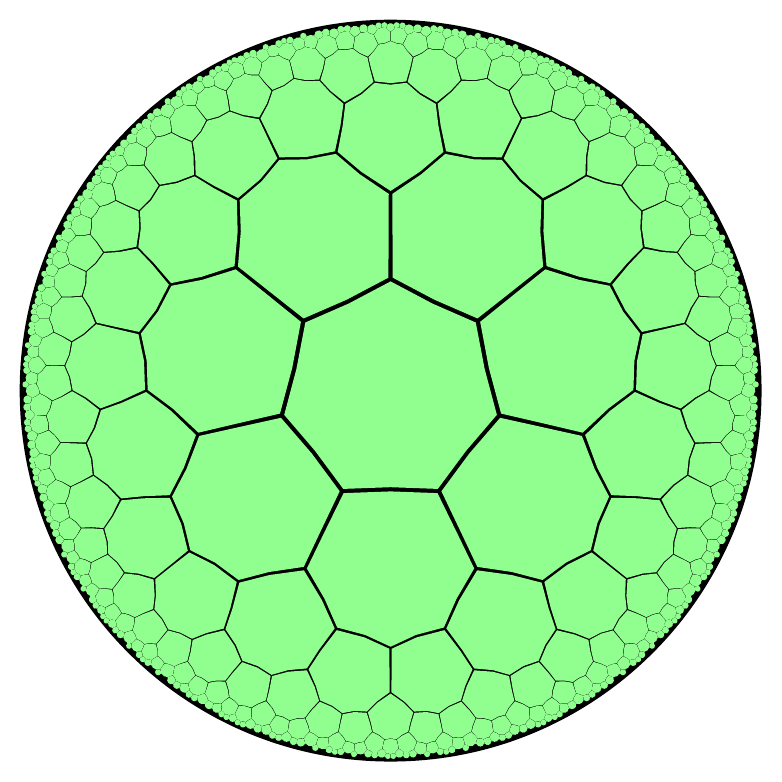}%
\includegraphics[width=.19\textwidth]{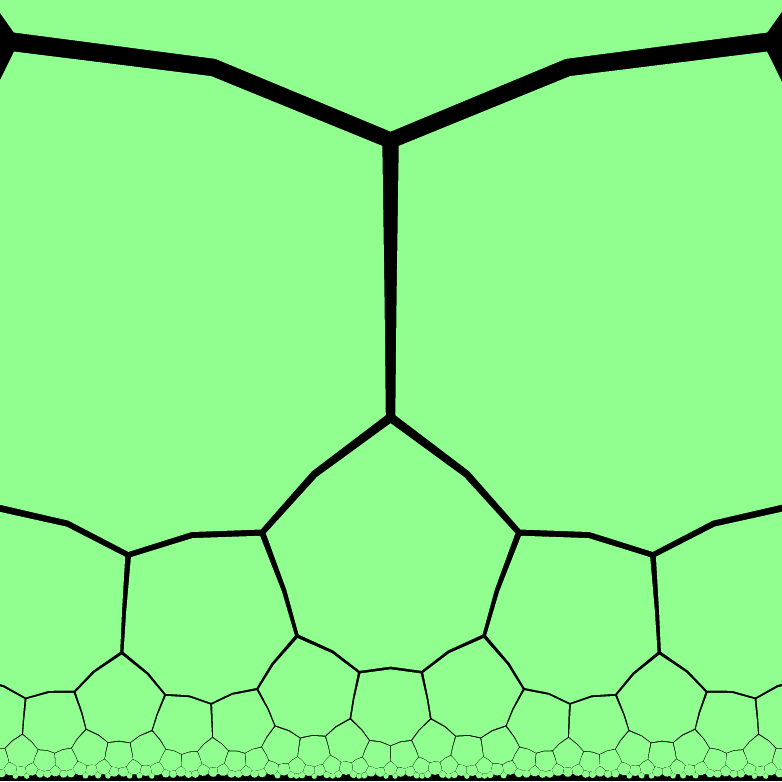}%
\includegraphics[width=.19\textwidth]{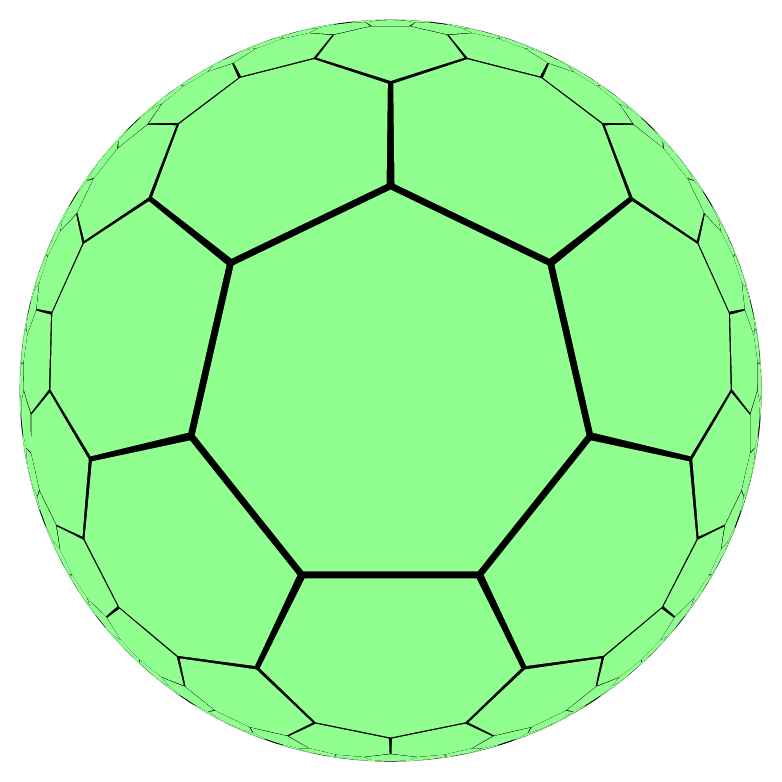}
\includegraphics[width=.19\textwidth]{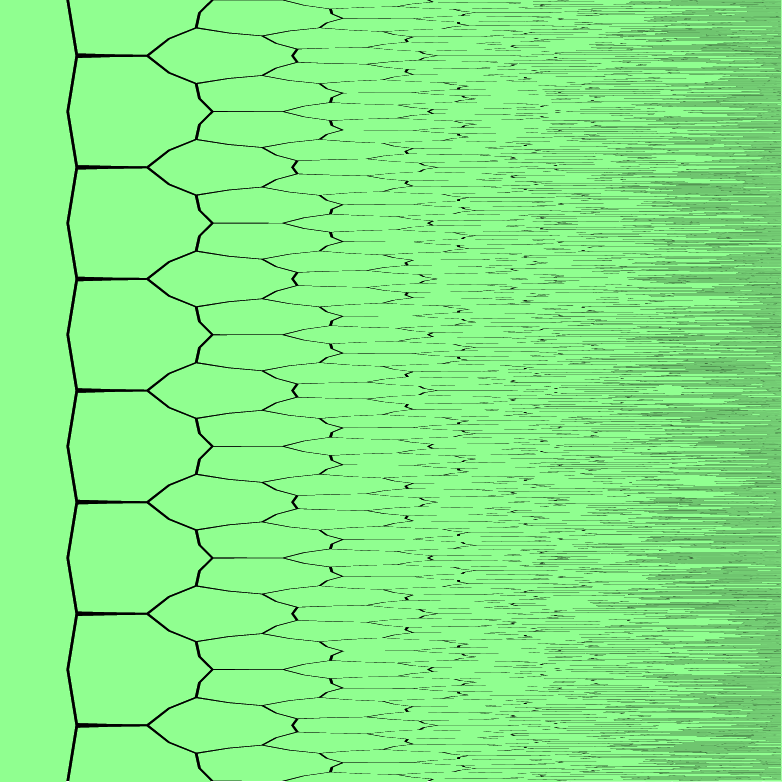}%
\end{center}
\caption{The \{7,3\} tessellation of the hyperbolic plane in the following models:\\
hyperboloid, Poincar\'e disk, upper half-plane, Beltrami-Klein disk, polar coordinates.\label{fig:tess}}
\end{figure}

\section{Hyperbolic geometry and representations}\label{sec2}

A geometry is defined by how points, lines, distances, angles, and isometries behave.
A \emph{representation} of a geometry is a method of representing points $p$ and isometries $f$. We need to
define and compute the following basic geometric objects and operations:

\begin{itemize}
\item {\sc Origin}, a constant $C_0$ representing a point which we consider the origin.
\item {\sc Translate} in the X direction by distance $x$ (returns an isometry $T^x$).
\item {\sc Rotate} around $C_0$ by an angle $\alpha$ (returns an isometry $R^\alpha$).
\item {\sc Apply} an isometry $f$ to a point $p$.
\item {\sc Compose} two isometries $f_1$ and $f_2$.
\item {\sc Invert} an isometry $f$.
\item {\sc Distance} of a given point $p$ to the origin $C_0$.
\end{itemize}

\paragraph{Spherical geometry}
The most straightforward non-Euclidean geometry is spherical geometry.
The $d$-dimensional sphere is $\bbS^d = \{x \in \bbR^{d+1}: g^+(x,x)=1\}$, where
$g$ is the inner product, $g^+((x_1,\ldots,x_{d+1}), (y_1,\ldots,x_{d+1})) = x_1y_1 + \ldots + x_{d+1}y_{d+1})$.
We consider $C_0=(0,\ldots,0,1) \in \bbS^d$ to be the origin of $\bbS^d$. 
Isometries of the sphere are exactly the isometries of $\bbR^{d+1}$ which map 0 to 0, i.e.,
orthogonal matrices in $\bbR^{(d+1) \times (d+1)}$. 
It is straightforward to compute the basic translations and rotations, for example,
an isometry $T^\alpha$ which moves the origin of $\bbS^2$ by $\alpha$ units in the direction of the first
coordinate can be written as $T^\alpha(x_1,x_2,x_3) = (x_1 \cos(\alpha) + x_3 \sin(\alpha), x_2, x_3 \cos(\alpha) - x_1 \sin\alpha)$.
The distance between the points $x,y \in \bbS^{d+1}$ is the length of the shortest arc
$\gamma \subseteq \bbS^{d+1}$ connecting $x$ and $y$, and can be computed using the formula $\acos\ g^+(x,y)$.

\paragraph{Minkowski hyperboloid (linear representation)}
The representation above lets a person with a basic knowledge of linear algebra work
with spherical geometry in an intuitive wall. We obtain hyperbolic geometry (in the Minkowski hyperboloid model) by applying
the same construction in pseudo-Euclidean space, using the Minkowski inner product
$g^-((x_1,\ldots,x_{d+1}), (y_1,\ldots,x_{d+1})) = x_1y_1 + \ldots + x_dy_d - x_{d+1}y_{d+1})$.
Many formulas of hyperbolic geometry are the same as the relevant formulas of spherical geometry, except that we
need to change the sign in some places (due to the change of sign in $g^-$), and use sinh and cosh instead of sin and cos when the argument represents distance.
The hyperbolic plane is $\bbH^d = \{x \in \bbR^{d+1}: x_{d+1}>0, g^-(x,x)=-1\}$. 
Again, we consider $C_0=(0,\ldots,0,1)$ to be the origin of $\bbH^d$. Isometries of $\bbH^d$
are the matrices $M \in \bbR^{(d+1) \times (d+1)}$ such that $g^-(Mx,My)=g^-(x,y)$, and 
$(MC_0)_{d+1}>0$. An isometry $T^\alpha$ which moves the origin of $\bbH^2$ by $\alpha$ units in the direction of the first
coordinate can be written as the Lorentz boost $T^\alpha(x_1,x_2,x_3) = (x_1 \cosh(\alpha) + x_3 \sinh(\alpha), x_2, x_3 \cosh(\alpha) + x_1 \sinh\alpha)$;
however, a rotation $R^\alpha$ of $\bbH^2$ by $\alpha$ around $C_0$ still uses cos and sin: 
$R^\alpha(x_1,x_2,x_3) = (x_1 \cos(\alpha) - x_2 \sin(\alpha), x_2 \cos(\alpha) + x_1\sin(\alpha), x_3)$.
The distance between the points $x,y \in \bbH^{d}$ is the length of the shortest arc
$\gamma \subseteq \bbH^{d}$ connecting $x$ and $y$ computed according to $g^-$, and can be computed using the formula $\acosh\ g^-(x,y)$.
For $x,y \in \bbH^{d}$, the midpoint $\midp(x,y)$ is the point in the middle of this arc. We have $\midp(x,y) = \frac{x+y}{-\sqrt{g^-(x+y, x+y)}}$.
We can use this model directly (i.e., represent points as their coordinates in $\bbH^d$ and isometries as linear transformation matrices), which we will call the
\emph{linear} representation.

\paragraph{Other models of $\bbH^2$, and polar representation.}
In the Minkowski hyperboloid model, every point has three coordinates, while two are sufficient.
Other common models (projections to $\bbR^2$) of $\bbH^2$ include:
\begin{itemize}
\item Beltrami-Klein disk model: the point $x=(x_1,x_2,x_3) \in \bbH^2$ is mapped to $P(x) = (\frac{x_1}{x_3}, \frac{x_2}{x_3})$.
This maps $\bbH^2$ to the inside of the unit disk in $\bbR^2$.
\item Poincar\'e disk model: the point $x=(x_1,x_2,x_3) \in \bbH^2$ is mapped to $K(x) = (\frac{x_1}{x_3+1}, \frac{x_2}{x_3+1})$.
Again, this maps $\bbH^2$ to the inside of the unit disk in $\bbR^2$.
\item Upper half-plane model: obtained from the Poincar\'e disk model by applying
a circle inversion which maps the inside of the Poincar\'e disk into the upper half-plane $\bbU^2$, which we interpret in terms
of complex numbers: $\bbU^2 = \{z \in \bbC: \Im(z) > 0\}$. The center point of $\bbU^2$ is $i \in \bbU^2$.
\item Native polar coordinates: in this model, we use coordinates $(\phi,r)\in\bbP$, where $r$ is the distance from $C_0$ and $\phi$
is the angle. In other words, the coordinates $(\phi,r)$ in native polar coordinates correspond to the point $R^\phi(T^r(C_0)) \in \bbH^2$.
\end{itemize}

These models have natural analogs in higher dimensions.
Figure \ref{fig:tess} shows the hyperbolic plane in the models above. 
The pictures show the tessellation of $\bbH^2$ by regular hyperbolic heptagons;
all heptagons are of the same size and shape, but these sizes and shapes had to be changed by the projection to $\bbE^2$ used.

Native polar coordinates are the first example of an alternative representation of the hyperbolic plane, 
that is popular in network science applications \cite{bogu_internet,friedrich2023computing}. 
It is somewhat analogous to using latitude and longitude in spherical geometry.
It is straightforward to compute the distance between $(\phi_1,r_1)\in\bbP$ and $(\phi_2,r_2)\in\bbP$.
Indeed, $(0,r_1)\in\bbP$ corresponds to $(\sinh(r_1),0,\cosh(r_1)) \in \bbH^2$, and $(\phi_2,r_2)\in\bbP$.
corresponds to $(\cos(\phi_2)\sinh(r_2),\sin(\phi_2)\sinh(r_2),\cosh(r_2)) \in \bbH^2$; thus,
the distance $d$ satisfies the hyperbolic cosine rule
\begin{equation}
\cosh(d) = \cosh(r_1)\cosh(r_2) + \cos(\phi_2) \sinh(r_2) \cosh(r_2). \label{cosinedist}
\end{equation}
In general, we need to replace $\phi_2$ with $\phi_1$. 
If $(\phi_1,r_1)$ and $(\phi_2,r_2)$ are close, the following formula is numerically better \cite{tobias}:
\begin{equation}\cosh(d) = \cosh(r_1-r_2) + (1-\cos(\phi))\sinh(r_1)\sinh(r_2) \label{precisedist}
\end{equation}
While less useful in network science applications, 
we also need to represent arbitrary isometries, but this is also straightforward: $(\phi,r,\psi)$ represents the isometry $R^\phi T^r R^\psi$.
In higher dimensions, we need to replace $\phi$ and $\psi$ with isometries of $\bbS^{d-1}$.

\def\clm{}
\def\cliff#1{Cl(#1)}
\def\cliffv{\cliff{\bbV}}
\def\cliffi#1#2{Cl^{#1}(#2)}
\def\cliffiv#1{\cliffi{#1}{\bbV}}
\def\cconj#1{\overline{#1}}

\paragraph{Clifford algebras, and mixed representation.}
While using $d \times d$ matrices is a straightforward method of representing rotations of $\bbR^d$ (and thus also isometries of $\bbS^{d-1}$),
it is often advantageous to use other representations. For orientation-preserving isometries of $\bbR^3$, quaternions are commonly used in computer graphics.
This section explains Clifford algebras that generalize this construction.

Let $\bbV$ be $\bbR^d$ with inner product $g$, and let $e_1, \ldots, e_d$ be the unit vectors of $\bbV$. The points of $\bbV$ can be written as $x_1e_1 + \ldots + x_de_d$.
The free algebra over $\bbV$,
$T(\bbV)$, is the vector space whose basis is the set of all sequences of $e_i$ (including the empty sequence, denoted by 1).
The elements of $T(\bbV)$ are added and multiplied in a natural way, for example:
\[(3+2e_2e_1) \clm (1+2e_1) = 3+2e_2e_1+6e_1+4e_2e_1e_1\]

Note that $\bbV$ is a subspace of $T(\bbV)$ and that this multiplication is associative but not commutative.
Addition is associative and commutative, and multiplication is distributive over addition.

The Clifford algebra $\cliffv$ is obtained from
$T(\bbV)$ by identifying elements according to the following rule: for $u, v \in \bbV$, $uv+vu = 2g(u,v)$ ($\star$). We perform all the
identifications that follow from this rule and the associativity/commutativity/distributivity rules. In particular, for $\bbV=(\bbR^{d+1},g^-)$
we have $e_ie_j=-e_je_i$ and $e_ie_i=1$ for $i \leq d$, and $e_ie_i=-1$ for $i=d+1$. Thus, $\cliffv$ is a $2^{d+1}$-dimensional space
(using the rules above, we can rewrite any element of $\cliffv$ using only products of $e_i$ which are ordered and have no repeats).

For $x \in \cliffv$, $\cconj{x}$, called the \emph{conjugate} of $x$, is defined as follows: $\cconj{1}=1$,
$\cconj{e_i}=-e_i$, $\cconj{x+y}=\cconj{x}+\cconj{y}$, $\cconj{x \clm y} = \cconj{y} \clm \cconj{x}$.

Let $v$ be a non-zero vector in $\bbV$. Any vector $w \in \bbV$ can then be decomposed as $w = av+u$, where $u$ is orthogonal to $v$.
We have
\( v\clm w\clm\cconj{v} = v\clm(av+u)\clm\cconj{v} = a(v\clm v\clm\cconj{v}) + v\clm u\clm\cconj{v} = 
g(v,v) a \cconj{v} - v \clm \cconj{v} \clm u = -g(v,v) a v + v \clm v \clm u = g(v,v) (u-av)\).
For $\bbV=(\bbR^{d+1},g^-)$, if $v,w \in \bbH^d$, the operation $w \mapsto v\clm w\clm\cconj{v}$
is exactly the point reflection of $\bbH^{d-1}$ in $v$; and if $g^-(v,v)=1$, it is the reflection in the hyperplane orthogonal to $v$.

For an $x \in \cliffv$, let us denote the operation $w \mapsto x\clm w\clm\cconj{x}$ by $M(v)$. It is easy to check that
$M(xy)$ is the composition of $M(x)$ and $M(y)$. Since the basic translations and rotations of $\bbH^d$ can be obtained as compositions
of two operations from the last paragraph, they can be represented as $M(x)$ for some $x \in \cliffv$; furthermore, every
orientation-preserving isometry of $\bbH^d$ can be obtained as a composition of basic translations and rotations, thus, also
$M(x)$ for some $x \in \cliffv$. It is easy to check that, since the number of basic reflections is even, we only use
$\cliffiv{[0]}$, which is the subspace of $\cliffv$ whose base is all the products of even number of $e_i$'s.

Therefore, for $\bbH^d$ (and $\bbS^d$), we have a representation of orientation-preserving isometries which uses only
$2^d$ real numbers. For $d \leq 5$ this is less than ${(d+1)}^2$ we would have to use in the $\bbR^{(d+1) \times (d+1)}$ representation.

\paragraph{Reduced representation.}
In the previous paragraph, we represented the points $x \in \bbH^d$ using the Minkowski hyperboloid model
and the isometries using $\cliffiv{[0]}$ (the \emph{mixed} approach).
Another possible representation is to represent $x \in \bbH^d$ as $y \in \cliffiv{[0]}$
such that $M(y)(C_0)=x$ and $M(y)(x')=C_0$, where $x'$ is the point such that $C_0=\midp(x,x')$
(i.e., this isometry moves $C_0$ to $x$ without introducing any rotation). This $y$ can be computed as $x''C_0$, where
$x''=\midp(C_0,x)$.

Reduced representation represents every point $x \in \bbH^d$ as $y_1e_1e_{d+1} + y_2 e_2e_{d+1} + \ldots + y_d e_d e_{d+1} + y_{d+1}$
(only $d+1$ out of $2^d$ coordinates are used). The coordinates $y=(y_1, y_2, \ldots y_{d+1})$ can be interpreted as coordinates
on the Minkowski hyperboloid of the point $x''$. Since $x''=\midp(C_0,x)$, the distance of $x''$ from $C_0$ is half the distance
of $x$ from $C_0$; this will be a serious numerical advantage. We call this representation \emph{reduced}.
Another important property is that, by projecting $y$ as in the Beltrami-Klein disk model, we obtain $x$ in the Poincar\'e disk model, $K(y)=K(x'')=P(x)$.

\paragraph{Half-plane representation.}
Another representation of $\bbH^2$ uses the upper half-plane $\bbU$ for points and matrices $A=\left(\begin{array}{cc}a & b \\ c & d\end{array}\right)$, where $a,b,c,d \in \bbR$,
for isometries.
We apply $A$ to $z \in \bbU^2$ as follows: ${\textsc{Apply}}(A,z) = (ax+b)/(cx+d)$. Note that, for $\alpha \in \bbR$ such that $\alpha \neq 0$, $A$ and $\alpha A$ represent the same isometry.
The \emph{normalized} one is the one which has determinant $ad-bc=1$. Therefore, 
the set of isometries corresponds to the set of 2$\times$2 matrices over reals with determinant 1, which is called the \emph{special linear group}
over $\bbR$, SL$(2,\bbR)$.

This representation of isometries is essentially equivalent to the Clifford algebra representation up to the base change. In particular, 
if the Clifford algebra representation of an isometry is $k_0 + k_2 e_1e_3 + k_1 e_2e_3 + k_3 e_1e_2$, then its SL$(2,\bbR)$ representation is
$\left(\begin{array}{cc}k_0-k_2 & k_1+k_3 \\ k_1-k_3 & k_0+k_2\end{array}\right)$.

\paragraph{Half-space representation.}
While half-plane representation works for $\bbH^2$, a similar representation exists for $\bbH^3$. We will be using quaternions $\bbH$\label{notation:quaternion}: a four-dimensional space
over reals, with the four basis vectors called 1, $i$, $j$, and $k$, multiplied according to rules $i^2=j^2=k^2=-1$, $ij=k$, $jk=i$, $ki=j$.
(To avoid confusing $\bbH$ with the hyperbolic space $\bbH^d$, note that the standard notation for quaternions, $\bbH$, has no index.)
The points are represented in the Poincar\'e half-space model, using quaternions $x \in \bbH$ such that the $j$-part of $x$ is positive, and the $k$-part of $x$ is 0.
The center point is $j \in \bbH$. The isometries are represented as SL($2,\bbC$), that is, 2$\times$2 matrices over $\bbC$ with determinant 1.
Applications are performed using the same formula: ${\textsc{Apply}}(A,x) = (ax+b)/(cx+d)$.
Again, this representation is essentially equivalent to the Clifford algebra: 
\(k_0+k_3e_1e_2+k_5e_1e_3+k_6e2e_3+k_9e_1e_4+k_{10}e_2e_4+k_{12}e_3e_4+k_{15}e_1e_2e_3e_4\) is equivalent
to \[\left(\begin{array}{cc}
k_0-k_9+k_{15}i-k_6i & k_3+k_{10}-k_5i-k_{12}i \\
k_{10}-k_3+k_{12}i-k_5i & k_0+k_9+k_6i+k_{15}i
\end{array}\right).\]

\section{Representation variants}
In Section \ref{sec2}, we introduced the basic representations we will compare in our
study (linear, mixed, reduced, half-plane/half-space, polar, and generalized polar).
There are also multiple methods of dealing with numerical errors. As an example, consider a
point $x=(x_1,\ldots,x_{d+1})$ in the Minkowski hyperboloid $\bbH^d$. We should have
$x_{d+1}^2=1+x_1^2+\ldots+x_d^2$. However, if we apply several of representation operations
to compute the point $x$, it may happen that this equation is not true as a result of numerical errors.
We consider the following variations:

\begin{itemize}
\item {\bf Invariant.} Do not do anything. Hope that numerical errors do not build up.

\item {\bf Careless.} Here we consider $\alpha x$, for any $\alpha \in \bbR$ other than 0,
to be a correct representation of $x \in \bbH^d$.

\item {\bf Flattened.} We normalize in another way: we multiply $x \in \bbH^d$ by $1/x_{d+1}$.
This lets us conserve memory since the $d+1$-th coordinate in our representation will always
equal 1. (This is effectively the Beltrami-Klein model.)

\item {\bf Forced.} Normalize the output after every computation.

\item {\bf Weakly Forced.} Try to normalize the output after every computation, but do not do it
if the norm could not be computed due to precision errors.

\item {\bf Binary.} In careless, values may easily explode and cause underflow/overflow; avoid this
  by making the leading coordinate in $[0.5, 2)$ range (by multiplying by powers of 2, which is
  presumably fast).
\end{itemize}

Similar variants can also be applied to the Clifford representation of points and isometries.
For example, flattened reduced representation is effectively the Poincar\'e disk model.
Furthermore:
\begin{itemize}
\item In linear representations, matrices can be \emph{fixed} by replacing them with correct orthogonal matrices close
  to the current computation. HyperRogue \cite{hyperrogue} uses this method;
  not applying such fixes would have a visible effect of the visualization becoming visibly stretched after
  the user moves sufficiently far away from the starting point (interestingly, it tends to
  fix itself when the user moves back towards the starting point). We call non-fixed representations
  \emph{linear-F}, and fixed representations \emph{linear+F}.
\item In \emph{polar1}, we always use the basic cosine rule (\ref{cosinedist}). In \emph{polar2},
we use a better formula when the angles are close or opposite to each other (\ref{precisedist}).
In 2D, we can use angles (\emph{angles} variant), but forcing angles into $[-\pi,\pi]$ may be needed to prevent explosion (\emph{mod} variant).
In general, we have a choice of representation for $\bbS^{d-1}$. We use 
\emph{forced} or \emph{invariant} reduced representations.
\item In the Clifford representation, the \emph{gyro} variant splits the isometries into 
  the translational part (which is flattened, making it equivalent to the Poincar\'e disk model)
  and the rotational part (for which 'invariant' is used). This fixes the problem 
  with full flattening where rotations by 180° are flattened to infinity. This is
  inspired by gyrovectors \cite{ungar} and is essentially doing the computation in the
  Poincar\'e disk model, a popular representation \cite{hyperbolica,unimodel,ccgcn}.
\end{itemize}

\begin{figure}
\begin{center}
\includegraphics[width=.14\textwidth]{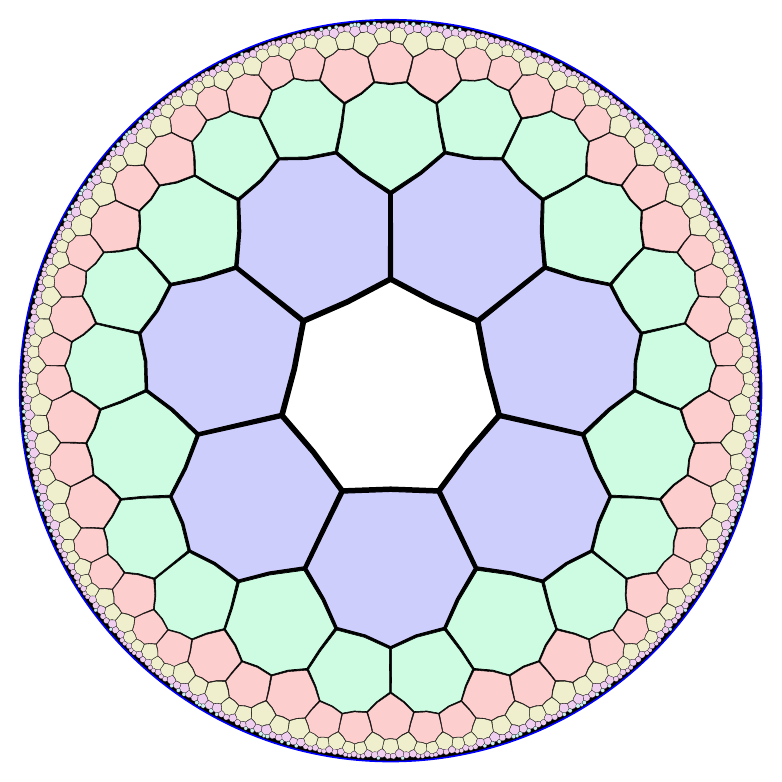}%
\includegraphics[width=.14\textwidth]{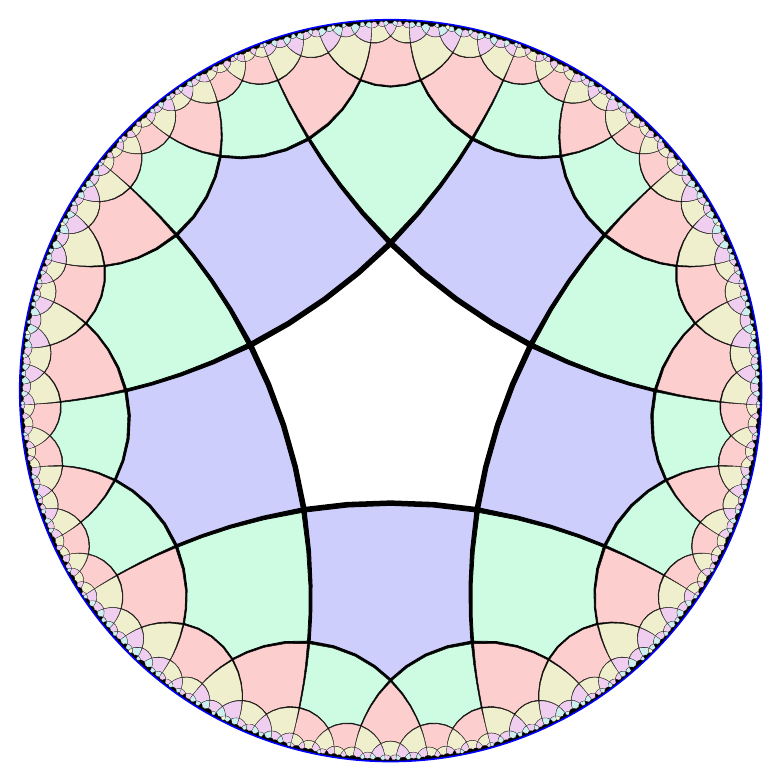}%
\includegraphics[width=.14\textwidth]{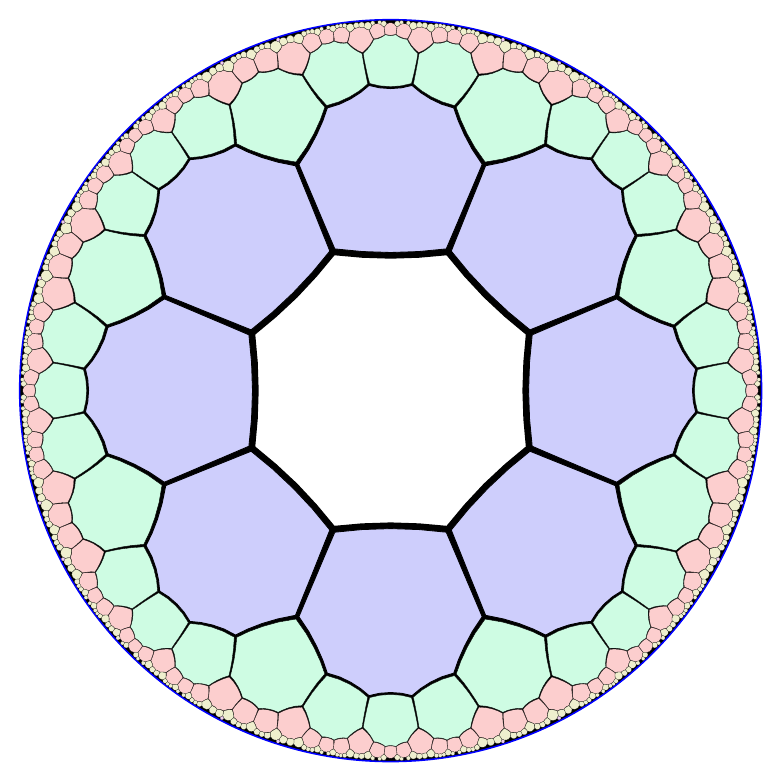}%
\includegraphics[width=.14\textwidth]{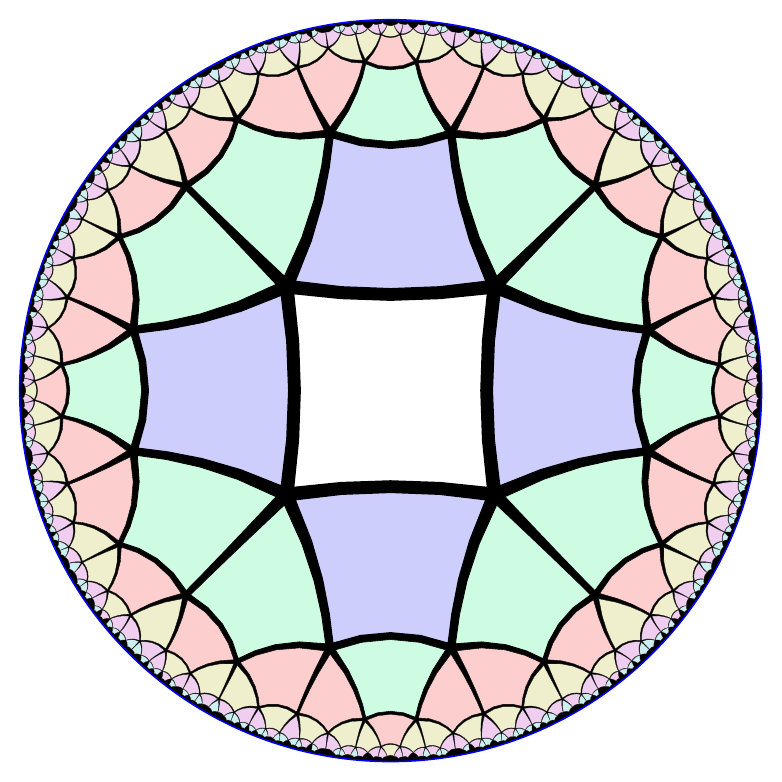}%
\includegraphics[width=.14\textwidth]{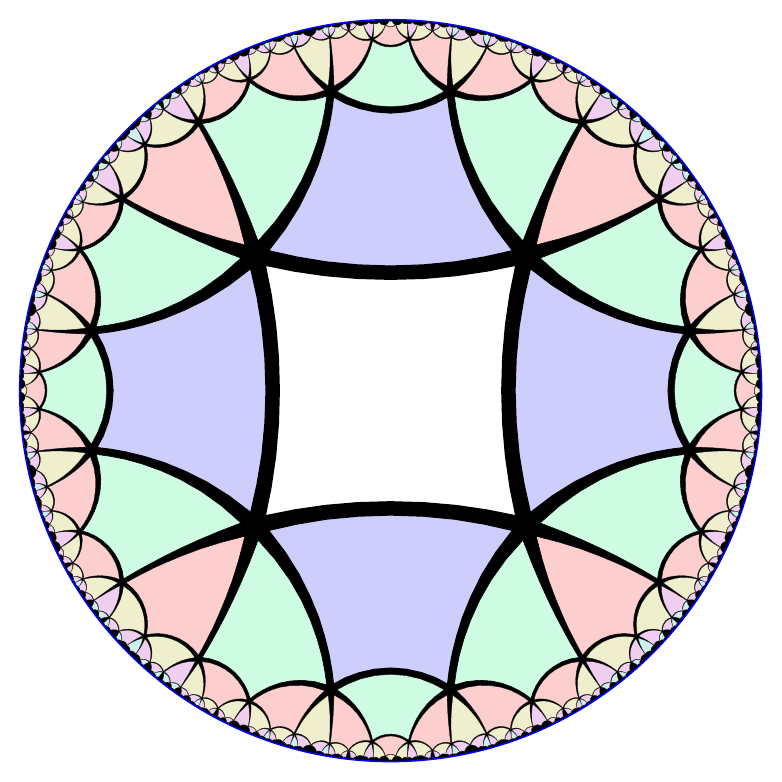}%
\includegraphics[width=.14\textwidth]{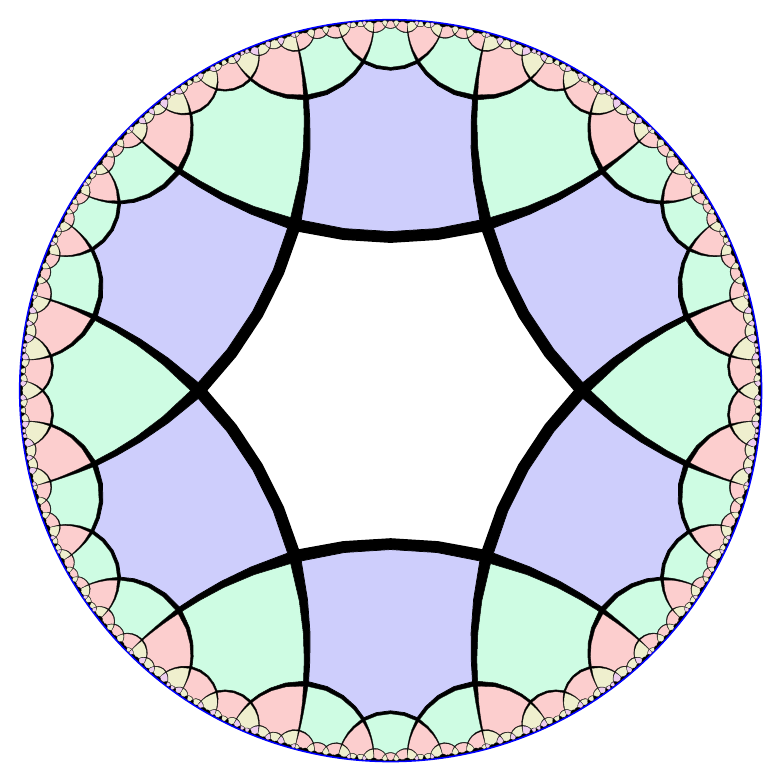}%
\includegraphics[width=.14\textwidth]{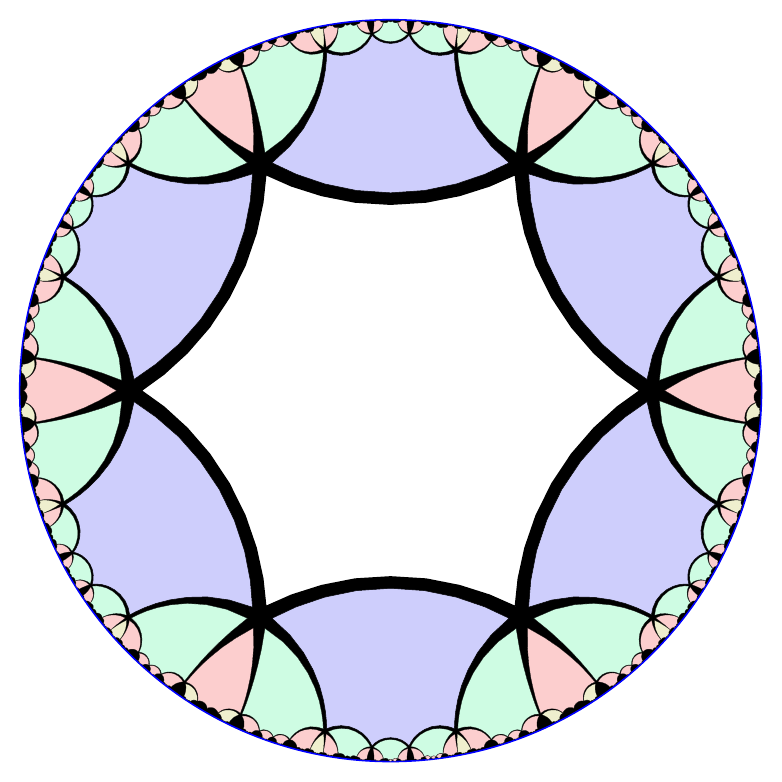}%
\end{center}
\caption{Tessellations of $\bbH^2$ used in our experiments.
From left to right: $\{7,3\}$, $\{5,4\}$, $\{8,3\}$, $\{4,5\}$, $\{4,6\}$, $\{6,4\}$, $\{6,6\}$.
\label{fig:vartess}}
\end{figure}

\section{Tessellations}

The hyperbolic plane can be tessellated with regular $p$-gons of the same size, such that $q$
of them meet in every vertex, for all $p,q$ such that $\frac{1}{p}+\frac{1}{q} < \frac{1}{2}$.
Such a tessellation has Schl\"afli symbol $\{p,q\}$. Figure \ref{fig:tess} shows the \{7,3\} tessellation.

Tessellations can be used to avoid numerical precision errors. With every tile $t$, we assign
an isometry $X_t$ that maps the central tile to $t$ (rotation is chosen arbitrarily).
If a point $p$ is in tile $t$, we can write $p = X_{t} p_0$, where $p_0$ is coordinates of $p$
relative to the center of $t$. We can now represent $p$ as $(t,p_0)$ -- since $p_0$ is
close to the origin, this avoids the problem of 
numerical precision issues being significantly higher for points
far away. A similar method can be used for isometries.

If two tiles $t$, $t'$ are adjacent, the isometry 
$X_{t,t'}=X_t X^{-1}_{t'}$ which maps $t'$-relative coordinates to $t$-relative coordinates 
equals $R^{\alpha} T^x R^{\beta}$, where the angles $\alpha$ and $\beta$ correspond to the chosen
orientations of $t$ and $t'$, and $x$ is the distance between two adjacent tiles, which can be
computed using hyperbolic trigonometry. The adjacency structure of tiles can be computed combinatorially
\cite{wpigroups,gentes}.

In our tests, we use tessellations either to produce tests with known correct answers or (as described above) to enhance the numerical precision. See Figure \ref{fig:vartess} for
the tessellations we use.

\section{Tests}

We compare representations on the following tests (some parameterized by $d$). 

\begin{description}
\item[LoopIso]
In this test, we construct a path $t_0,\ldots,t_k$ in the tiling by always moving to a random
adjacent tile until we get to a tile $d$ afar; then, we return to the start (also randomly,
may stray further from the path). We compose all the relative tile isometries 
$X_{t_k, t_{k-1}} \ldots (X_{t_2, t_1} X_{t_1,t_0})$ into $f$, which theoretically should 
equal identity. We see if $f(C_0) = C_0$. The test result is the first
distance $d$ for which $f(C_0)$ is not found to equal $C_0$ (distance $> 0.1$).


\item[LoopPoint]
Same as LoopIso, but we apply the consecutive isometries to point right away.
We compute $h = X_{t_k, t_{k-1}} \ldots (X_{t_2, t_1} (X_{t_1,t_0} C_0))$. We see if $h=C_0$.


\item[AngleDist]
We construct a random path $t_0,\ldots,t_d$. This time, we do not loop.
We compute $h = X_{t_d, t_{d-1}} \ldots (X_{t_2, t_1} (X_{t_1,t_0} C_0))$.
We check whether the angle and distance of $h$ from $C_0$ have been computed correctly.
In the variant \emph{AngleDist2}, we multiply the matrices in the opposite order.
The correct angle and distance are computed using high-precision floating point numbers.
The test result is the first distance $d$ for which the error exceeds 0.1.


\item[Distance]
We compute the distance between two points in distance $d$ from the starting point.
Such computations are of importance in social network analysis applications.
The angle between them is very small (similarity), or close to $180^\circ$ (dissimilarity),
close to $1^\circ$ (other). The test result is the first distance $d$ for which the error exceeds 0.1.


\item[Walk]
This test is based on an effect of numerical precision issues that is most visible in HyperRogue \cite{hyperrogue}.
After walking in a small line, it can often be clearly observed that we
have ``deviated'' from the original straight line. This test checks how long we can walk until this happens.

We construct an isometry $A$ representing a random direction. In each step, we compose this isometry with a translation ($A := AT^{1/16}$).
Whenever the point $AC_0$ is closer to the center of another tile, we rebase to that new tile.

For a test, we do this in parallel with two isometries $A$ and $B$, where $B$ = $AT^{1/32}$. We count the number of steps
until the paths diverge. In the \emph{WalkGood} variant, we instead compare to the result obtained
using high-precision floating point numbers.


\item[Close]
Here, we see whether minor errors accumulate when moving close to the center. Like in loop tests,
we move randomly until we reach distance $d+1$,
after which we return to the start (always reducing the distance). After each return to the start, we check if the representation is
still fine; if yes, we repeat the loop, letting the errors accumulate over all loops. We stop when the error is high or after 10000 steps.
The variants \emph{Close} and \emph{CloseInverse} differ in the order of multiplying $X$ matrices.



\end{description}

\section{Experimental Results}

\begin{figure}[h!]
\begin{center}
\includegraphics[width=\textwidth]{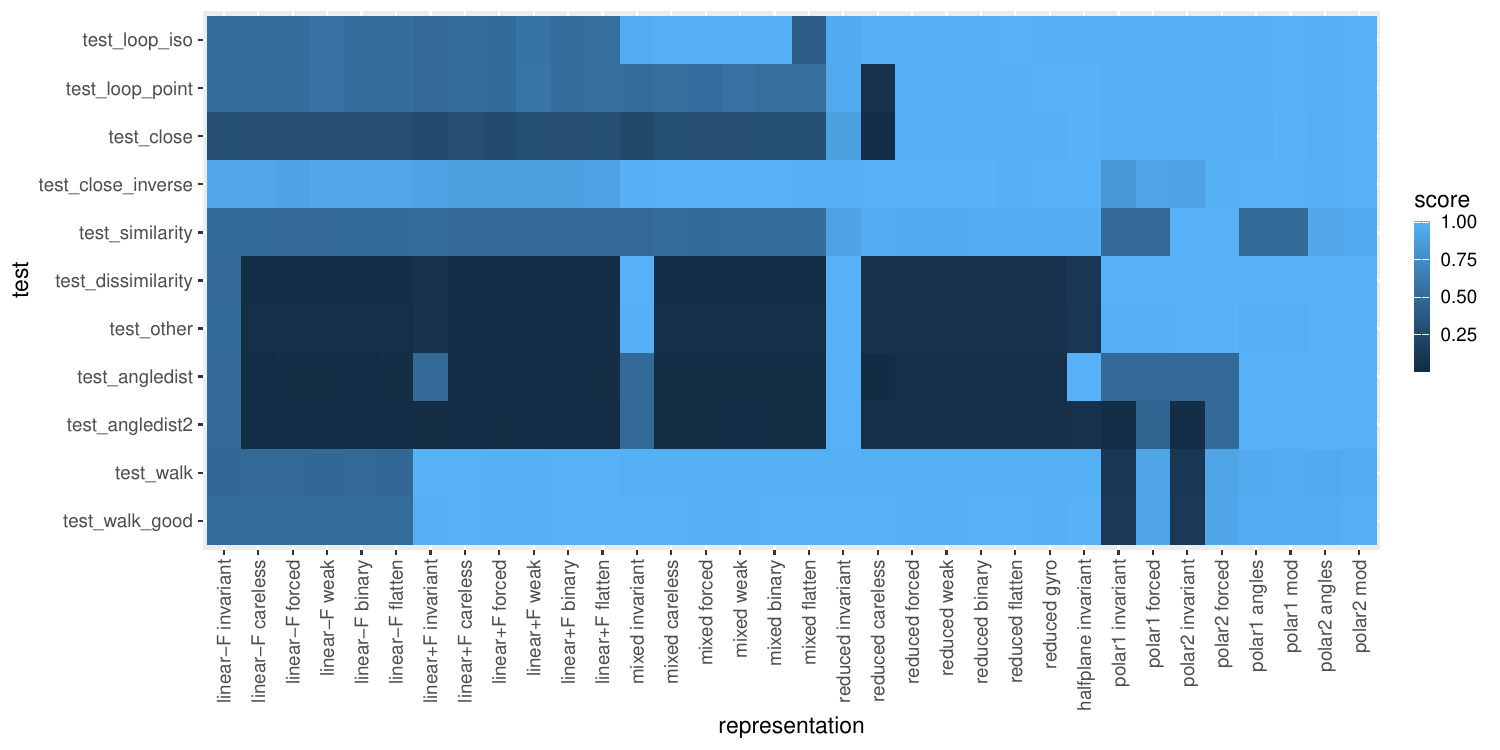}
\includegraphics[width=\textwidth]{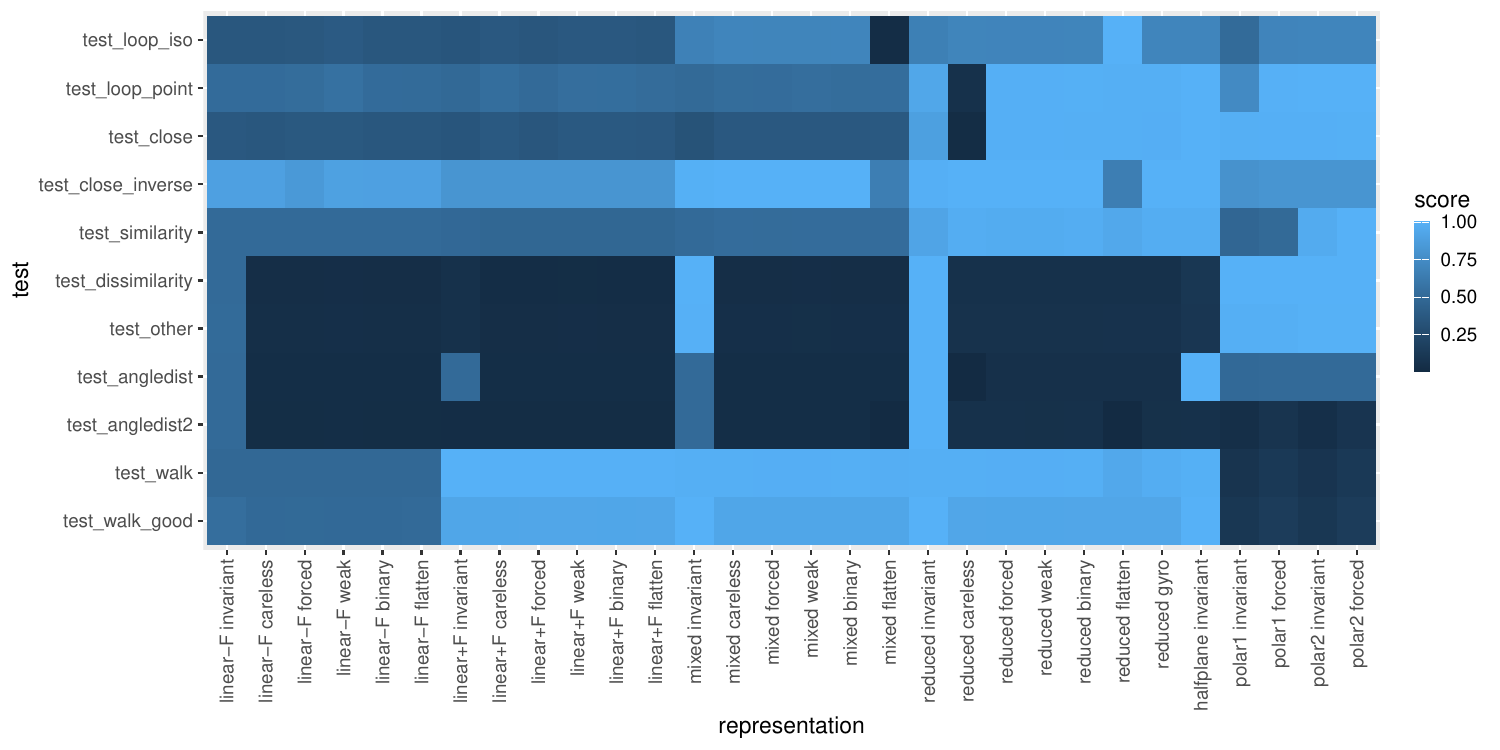}%
\end{center}
\caption{Aggregate results of our tests, in 2 dimensions (above) and 3 dimensions (below).\label{fig:aggregate}
For each test, the score of the best representation is normalized to 1. 
}
\end{figure}

For our implementation and experimental results, see \cite{numerical_fig}.
Our implementation is written in C++. All computations use the IEEE754 double precision format.
The GMP library is used for the high-precision comparisons.

Figure \ref{fig:aggregate} presents the heatmap of averages of running each representation
on each test 20000 times. We observe significant differences between the performance of
various representations. Some of these experimental results can be quite easily explained.
It is clear from Figure \ref{fig:tess} that the Beltrami-Klein disk model 
is not good
at representing points far away from $C_0$.
In the Klein disk model, a point in distance $d$ from $C_0$ is mapped to a point
in distance $\tanh(d)$ from the center of the disk, which is $1-\Theta(\exp(2d))$.
Floating point numbers cannot express such a slight difference from 1. Conversely, in the Poincar\'e disk model, this is $1-\Theta(\exp(d))$.
Therefore, we can
expect the effective distance represented accurately in flattened reduced
to be double that of flattened linear.
This issue carries over to most computations in non-flattened representations, although
not all of them. This effect has been studied in \cite{achilles}.

In the invariant linear representation, the point $R^\alpha T^xC_0$ is represented as $(\cos(\alpha) \sinh(x), \sin(\alpha) \sinh(x), \cosh(x))$,
and since floating-point numbers are good at representing large numbers, we can recover $\alpha$ and $x$
even if $x$ is very large (as long as $\cosh(x)$ fits in the range of our floating-point type).
This makes invariant linear significantly better than the Poincar\'e disk model (flattened reduced) in some experiments,
such as AngleDist. This relies on multiplying the matrices in the correct order and on the fact that we care only
about distances -- numerically computing $AB^{-1}$ for two isometries moving $C_0$ to two points $a$ and $b$ which
are closed to each other but far from $C_0$ is not likely to yield meaningful results. While
polar representations can represent even larger distances $x$, this does not carry over to computations
in our implementation, which have to compute $\exp(x)$, $\cosh(x)$ or $\sinh(x)$ anyway.

In three dimensions, flattened reduced fails due to some isometries of the $\{4,3,5\}$ honeycomb
not being representable (due to not having the unit component). The gyro variant fixes this issue.
In some tests, the lack of normalization causes the coordinates to quickly blow up exponentially, which
is avoided in normalized variants. This happens, e.g., in the \emph{Walk} test for polar invariant
(avoided in polar forced) and in the \emph{LoopPoint} test for reduced careless.

The invariant reduced representation shares both advantages of invariant and reduced representations.
In the half-plane and half-space models, we no longer have the problem of faraway points getting 
close to 1 (typically, they are complex numbers with the imaginary coordinate close to 0 instead).
So, these representations are quite good at representing large distances.
We can expect these representations to be highly accurate; our experiments support that. Polar representations are highly accurate, too.

\begin{figure}
\begin{center}
\includegraphics[width=0.9\textwidth]{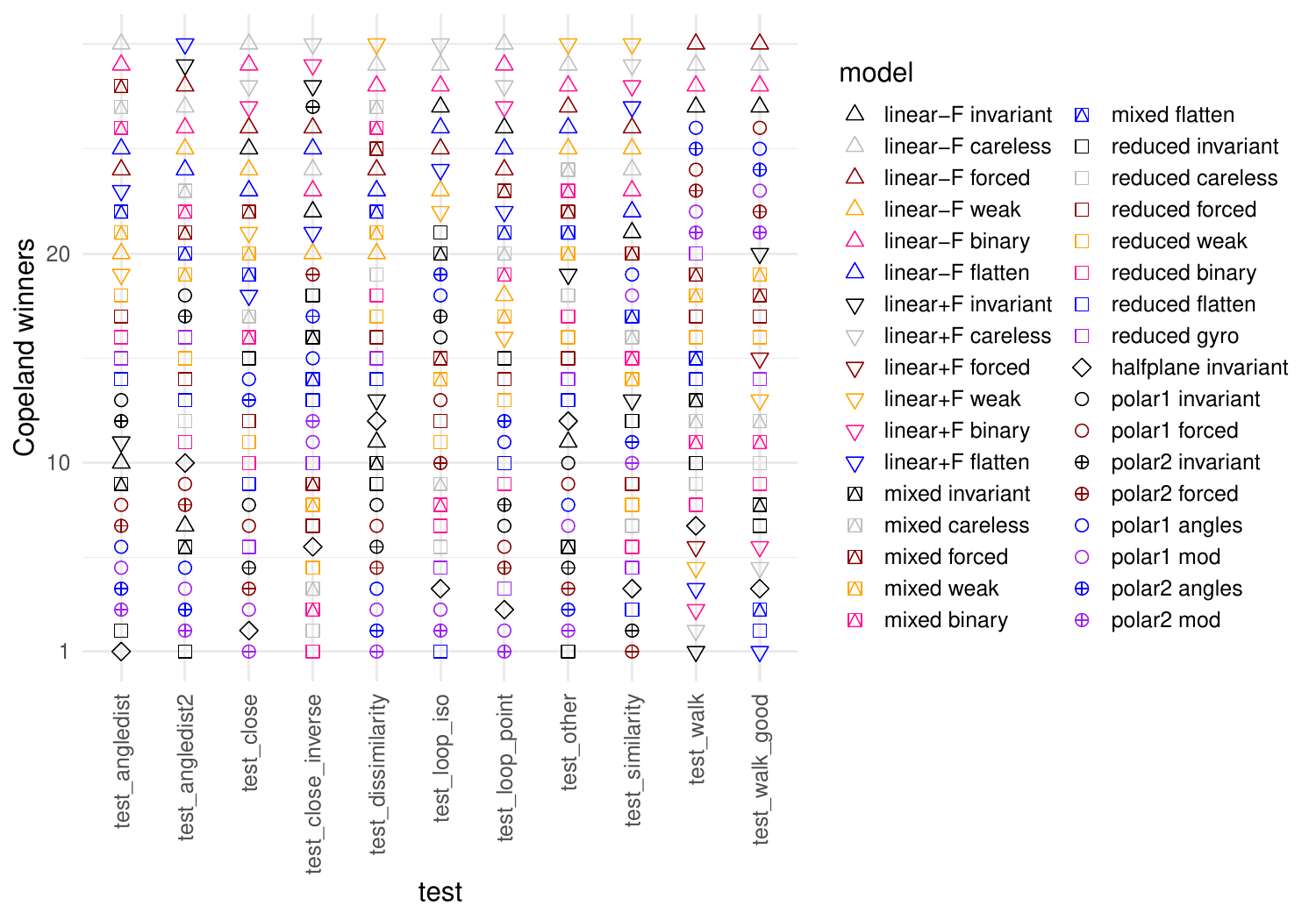}%
\end{center}
\caption{Condorcet rankings, in 2 dimensions. Better representations are closer to the bottom.\label{fig:condorcet2}}
\end{figure}

\begin{figure}
\begin{center}
\includegraphics[width=0.9\textwidth]{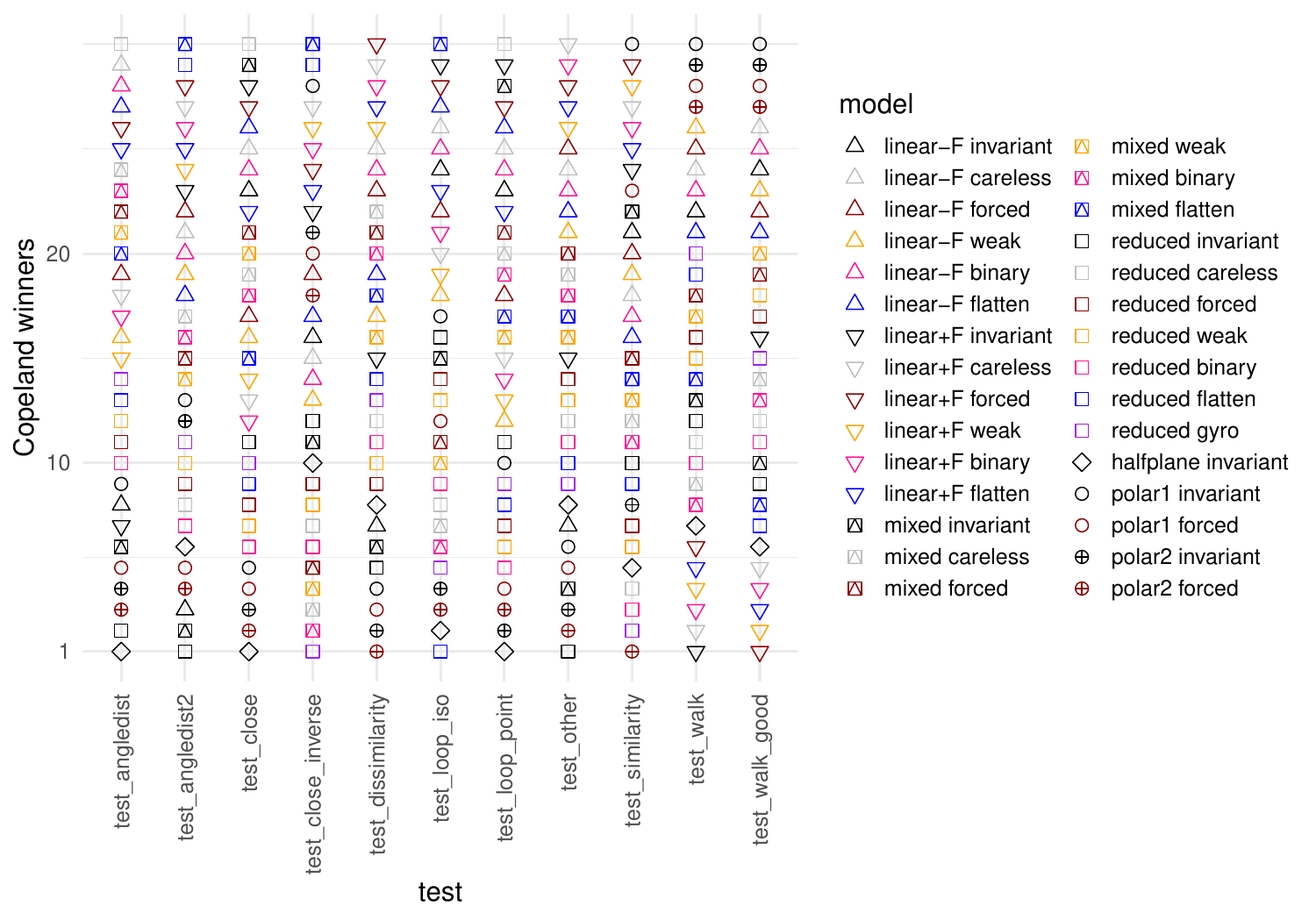}%
\end{center}
\caption{Condorcet rankings, in 3 dimensions.\label{fig:condorcet3}}
\end{figure}

Figures \ref{fig:condorcet2} and \ref{fig:condorcet3} depict the accurate rankings.
A representation $A$ wins against the representation $B$ if the probability that a randomly chosen simulation result obtained by A is greater
than a randomly chosen simulation result obtained by B exceeds 0.5. If that probability is equal to 0.5, we have a tie between A and B;
otherwise, A loses against B. We use the Condorcet voting rule (breaking ties using the Copeland rule \cite{copeland}) to obtain the ranking.
To compute the score for a given representation, we add 1 for every winning scenario, 0 for every tie, and -1 for every losing scenario.

The simulations are conducted under specific conditions, including 2000 iterations and the consideration of both 2D and 3D tessellations. These conditions are chosen to provide a comprehensive and representative view of the performance of each representation.

In many tests, polar2 mod representation wins, although the halfplane invariant is also very successful and
linear+F representations (especially invariant) perform well in walk tests.

\begin{figure}
\begin{center}
\includegraphics[width=\textwidth]{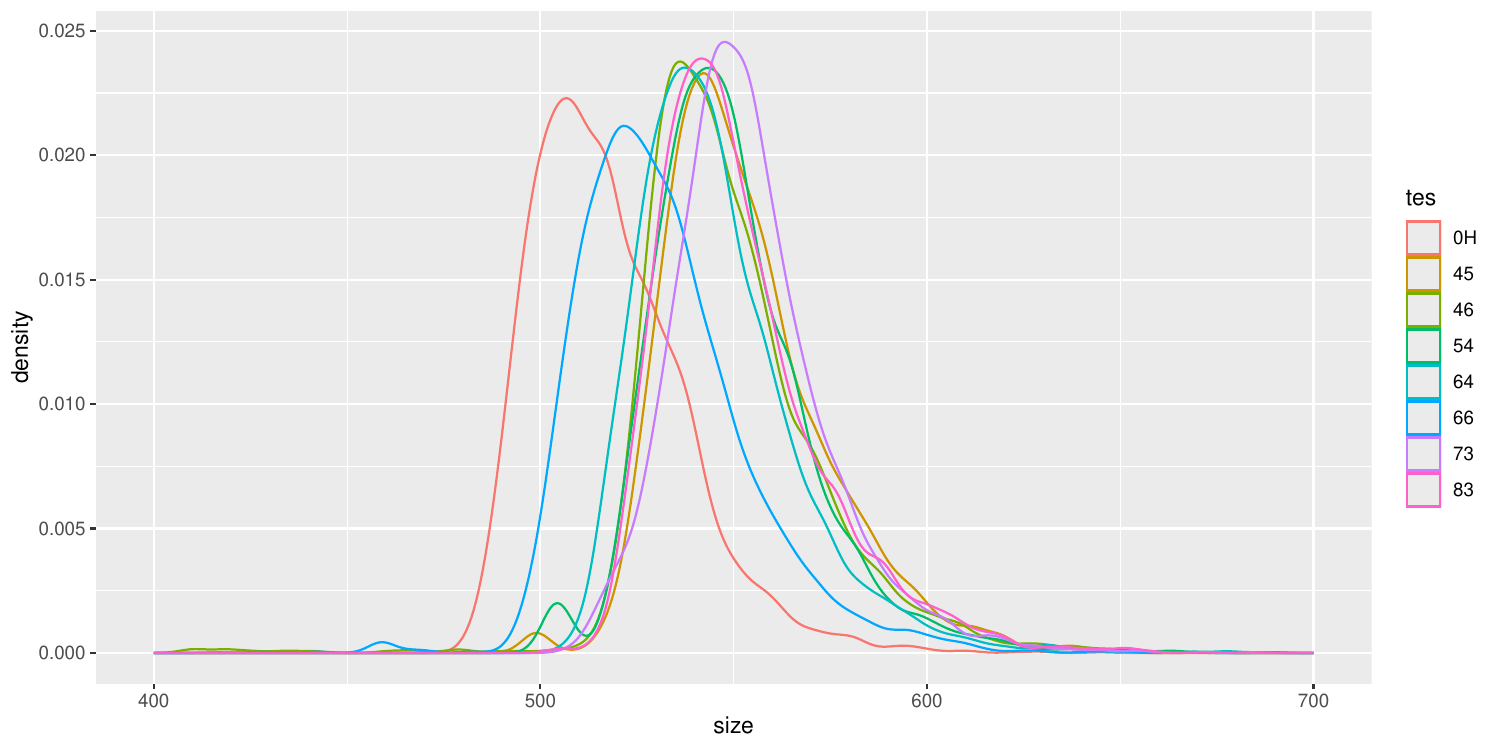}%
\end{center}
\caption{The result of the \emph{Walk} test, divided by tessellation. \label{fig:walkbig}}
\end{figure}

Figure \ref{fig:walkbig} depicts the density graph of
the results of the \emph{Walk} test, divided by tessellation.
We use the linear+F invariant representation since it achieves
the best results. For comparison, we also present the result of a
similar test without using any tessellation; this result uses 
the half-plane representation (linear+F invariant is half as good in this case).
From the tessellations we have considered, $\{7,3\}$ yields the best results;
other tessellations get slightly worse results except $\{6,6\}$ which is
significantly worse, and no tessellation is even worse.

\section{Comparison based on non-numerical advantages}
We focused on comparing various representations of hyperbolic geometry
concerning numerical precision. However, computational issues are not the only ones that should impact the final choice of the representation. That is why this section discusses the non-numerical advantages of particular representations.

The Poincar\'e model is often used to teach hyperbolic geometry, so many
people use it for this reason -- it is the one they know the best.
However, the linear representation is more natural for those who understand Minkowski geometry. For example, the formulas for
$T^x$ and $\midp(x,y)$ are straightforward in Minkowski hyperboloid,
but the respective Poincar\'e model formulas are not. Polar coordinates,
based on the idea of angle and distance from the origin, are
also intuitive for many people \cite{dziadICCS}. While the half-space models
have great numerical properties, their lack of symmetry makes them less intuitive.

One advantage of matrix representations is that
they can represent both orientation-preserving and orientation-reversing
isometries. To represent orientation-reversing isometries in Clifford algebras,
half-plane/half-space, or polar representations, we need to add an extra
bit of information -- whether the isometry should be composed with a mirror image.
This introduces some extra complexity.

In some applications (e.g., to work with non-compact honeycombs in $\bbH^3$ \cite{hhoney}) we need to represent not only
the \emph{material} points discussed so far, but also \emph{ideal} and \emph{ultra-ideal} points.
Intuitively, ideal points are the points on the boundary of the Beltrami-Klein model, 
and ultra-ideal points are outside of the boundary. In the Poincar\'e disk and
half-plane model, we can represent ideal points (they are on the boundary), but
we cannot represent ultra-ideal points (the points outside of the boundary are
better interpreted as alternative representations of material points).
Thus, only non-normalized (careless) linear/mixed can represent ultra-ideal points, 
while only non-normalized/flattened linear/mixed/reduced/half-plane/half-space can 
represent ideal points.

In some applications, it is useful to have a single implementation of 
all isotropic geometries (hyperbolic, Euclidean, and spherical, with varying curvatures). In \cite{unimodel,ccgcn}, the Poincar\'e disk model
(and its spherical analog, stereographic projection) is chosen for this reason.
However, polar coordinates also work in all isotropic geometries, and a similar general model can be obtained for
the hyperboloid ($\bbH^d$), sphere ($\bbS^d$), and plane model; furthermore, linear representations
can also be used for other Thurston geometries \cite{rtvizfinal}. From the representations we studied, it seems that only half-plane and half-space do not easily generalize beyond hyperbolic geometry.

There are applications of hyperbolic geometry where we focus on the neighborhood of a straight line $L$
rather than the origin point $C_0$. In such applications, in $\bbH^2$, it is useful to use coordinates $(x,y)$
where $x$ is the distance along $L$, and $y$ corresponds to the distance from $L$. Here,
$y$ can be simply the distance from $L$ (\emph{Lobachevsky coordinates}, where $T^x R^{\pi/2} T^y C_0$
gets coordinates $(x,y)$, analogous to longitude and latitude in spherical geometry), 
or in some other way, like Bulatov's conformal band model \cite{bulatovband,hyperboliccat}. 
Such a situation occurs in \cite{hyperboliccat}, where the conformal band model is used.
There are also potential applications in video games (racing along a straight line) and
data analysis (bi-polar data). This is a very specific potential application, with
major advantages in specific circumstances and no advantages otherwise, so we do not compare
such representations.

\section{Conclusions}
This paper aimed to determine which representation of hyperbolic geometry is best concerning numerical issues.
To this end, we compared five representations (linear, mixed, reduced, halfplane/halfspace, and generalized polar),
controlling for six variants of dealing with numerical errors (invariant, careless, flattened, forced, weakly forced, binary).
We conducted six tests capturing different scenarios leading to accumulating numerical imprecisions. Our results suggest that
polar representation is the best in many cases, although the halfplane invariant is also very successful.
It is known that numerical errors can be successfully combated by combining representation with tessellations,
so our research took that into account.
From the tessellations we have considered, $\{7,3\}$ yields the best results.
Fixed linear representations (especially invariant)
perform well in game-design-related scenarios (walk tests: how long we can walk until we observe that we have ``deviated''
from the original straight line).

\paragraph{Acknowledgments}
This work has been supported by the National Science Centre, Poland, grant UMO-2019/35/B/ST6/04456.

\def\ext#1{#1}
\bibliographystyle{plain}
\bibliography{master1}

\begin{thebibliography}{10}
\providecommand{\url}[1]{\texttt{#1}}
\providecommand{\urlprefix}{URL }
\providecommand{\doi}[1]{https://doi.org/#1}

\bibitem{ccgcn}
Bachmann, G., B\'{e}cigneul, G., Ganea, O.E.: Constant curvature graph
  convolutional networks. In: Proceedings of the 37th International Conference
  on Machine Learning. ICML'20, JMLR.org (2020)

\bibitem{tobias_alenex}
Bl\"asius, T., Friedrich, T., Katzmann, M., Krohmer, A.: Hyperbolic embeddings
  for near-optimal greedy routing. In: Algorithm Engineering and Experiments
  (ALENEX). pp. 199--208 (2018)

\bibitem{tobias}
Bl\"asius, T., Friedrich, T., Krohmer, A., Laue, S.: Efficient embedding of
  scale-free graphs in the hyperbolic plane. In: European Symposium on
  Algorithms (ESA). pp. 16:1--16:18 (2016)

\bibitem{unimodel}
Block, A., Skopek, O., Bachmann, G., Ganea, O., Bécigneul, G.: A universal
  model for hyperbolic, euclidean and spherical geometries.
  \url{https://andbloch.github.io/K-Stereographic-Model/} (2019)

\bibitem{bogu_internet}
Boguñá, M., Papadopoulos, F., Krioukov, D.: Sustaining the internet with
  hyperbolic mapping. Nature Communications  \textbf{1}(6),  1–8 (Sep 2010).
  \doi{10.1038/ncomms1063}

\bibitem{bulatovband}
Bulatov, V.: Conformal models of hyperbolic geometry.
  \url{http://bulatov.org/math/1003/band-stretch1.html\#(1)} (2010)

\bibitem{dhrgex}
Celi\'nska-Kopczy\'nska, D., Kopczy\'nski, E.: Discrete hyperbolic random graph
  model (2021)

\bibitem{gentes}
Celi\'nska-Kopczy\'nska, D., Kopczy\'nski, E.: Generating tree structures for
  hyperbolic tessellations (2021)

\bibitem{numerical_fig}
Celi\'nska-Kopczy\'nska, D., Kopczy\'nski, E.: Numerical aspects of hyperbolic
  geometry (2024),
  \url{https://figshare.com/articles/software/Numerical_Aspects_of_Hyperbolic_Geometry/25325338}
  (accessed Apr 13, 2024)

\bibitem{hyperbolica}
CodeParade: Hyperbolica (2022)

\bibitem{wpigroups}
Epstein, D.B.A., Paterson, M.S., Cannon, J.W., Holt, D.F., Levy, S.V.,
  Thurston, W.P.: Word Processing in Groups. A. K. Peters, Ltd., USA (1992)

\bibitem{achilles}
Floyd, W.J., Weber, B., Weeks, J.R.: The achilles' heel of o(3, 1)? Exp. Math.
  \textbf{11}(1),  91--97 (2002). \doi{10.1080/10586458.2002.10504472}

\bibitem{friedrich2023computing}
Friedrich, T., Katzmann, M., Schiller, L.: Computing voronoi diagrams in the
  polar-coordinate model of the hyperbolic plane (2023)

\bibitem{hyperbolicvr}
Hart, V., Hawksley, A., Matsumoto, E.A., Segerman, H.: Non-euclidean virtual
  reality {I}: explorations of $\mathbb{H}^3$. In: Proceedings of Bridges:
  Mathematics, Music, Art, Architecture, Culture. pp. 33--40. Tessellations
  Publishing, Phoenix, Arizona (2017)

\bibitem{hyperrogue}
Kopczy\'{n}ski, E., Celi\'{n}ska, D., \v{C}trn\'{a}ct, M.: Hyper{R}ogue:
  Playing with hyperbolic geometry. In: Proceedings of Bridges \ext{:
  Mathematics, Art, Music, Architecture, Education, Culture}. pp. 9--16.
  Tessellations Publishing, Phoenix, Arizona (2017)

\bibitem{hyperboliccat}
Kopczy{\'{n}}ski, E., Celi{\'{n}}ska-Kopczy{\'{n}}ska, D.: Conformal mappings
  of the hyperbolic plane to arbitrary shapes. In: Goldstine, S., McKenna, D.,
  Fenyvesi, K. (eds.) Proceedings of Bridges 2019: Mathematics, Art, Music,
  Architecture, Education, Culture. pp. 91--98. Tessellations Publishing,
  Phoenix, Arizona (2019), available online at
  \url{http://archive.bridgesmathart.org/2019/bridges2019-91.pdf}

\bibitem{rtvizfinal}
Kopczyński, E., Celińska-Kopczyńska, D.: Real-time visualization in
  anisotropic geometries. Experimental Mathematics  \textbf{0}(0),  1--20
  (2022). \doi{10.1080/10586458.2022.2050324}

\bibitem{lampingrao}
Lamping, J., Rao, R., Pirolli, P.: A focus+context technique based on
  hyperbolic geometry for visualizing large hierarchies. In: Proceedings of the
  SIGCHI Conference on Human Factors in Computing Systems. pp. 401--408. CHI
  '95, ACM Press/Addison-Wesley Publishing Co., New York, NY, USA (1995).
  \doi{10.1145/223904.223956}

\bibitem{copeland}
Maskin, E., Dasgupta, P.: The fairest vote of all. Scientific American
  \textbf{290}(3),  64--69 (2004)

\bibitem{munzner}
Munzner, T.: Exploring large graphs in 3d hyperbolic space. {IEEE} Computer
  Graphics and Applications  \textbf{18}(4),  18--23 (1998).
  \doi{10.1109/38.689657}

\bibitem{hhoney}
Nelson, R., Segerman, H.: Visualizing hyperbolic honeycombs. Journal of
  Mathematics and the Arts  \textbf{11}(1),  4--39 (jan 2017).
  \doi{10.1080/17513472.2016.1263789}

\bibitem{nickel}
Nickel, M., Kiela, D.: Poincar\'{e} embeddings for learning hierarchical
  representations. In: Guyon, I., Luxburg, U.V., Bengio, S., Wallach, H.,
  Fergus, R., Vishwanathan, S., Garnett, R. (eds.) Advances in Neural
  Information Processing Systems 30, pp. 6341--6350. Curran Associates, Inc.
  (2017),
  \url{http://papers.nips.cc/paper/7213-poincare-embeddings-for-learning-hierarchical-representations.pdf}

\bibitem{dziadICCS}
Osudin, D., Child, C., He, Y.H.: Rendering non-euclidean space in real-time
  using spherical and hyperbolic trigonometry. In: Rodrigues, J.M.F., Cardoso,
  P.J.S., Monteiro, J., Lam, R., Krzhizhanovskaya, V.V., Lees, M.H., Dongarra,
  J.J., Sloot, P.M. (eds.) Computational Science -- ICCS 2019. pp. 543--550.
  Springer International Publishing, Cham (2019)

\bibitem{reptradeoff}
Sala, F., De~Sa, C., Gu, A., Re, C.: Representation tradeoffs for hyperbolic
  embeddings. In: Proc.\ ICML. pp. 4460--4469. PMLR, Stockholmsm{\"a}ssan,
  Stockholm Sweden (2018), \url{http://proceedings.mlr.press/v80/sala18a.html}

\bibitem{ungar}
Ungar, A.: Gyrovector spaces and their differential geometry. Nonlinear
  Functional Analysis and Applications  \textbf{10} (01 2005)

\bibitem{mltiles}
Yu, T., De~Sa, C.M.: Numerically accurate hyperbolic embeddings using
  tiling-based models. In: Wallach, H., Larochelle, H., Beygelzimer, A.,
  d\textquotesingle Alch\'{e}-Buc, F., Fox, E., Garnett, R. (eds.) Advances in
  Neural Information Processing Systems. vol.~32. Curran Associates, Inc.
  (2019),
  \url{https://proceedings.neurips.cc/paper/2019/file/82c2559140b95ccda9c6ca4a8b981f1e-Paper.pdf}

\end{thebibliography}

\end{document}